\documentclass[aop,seceqn,citesort,MSNbibl,dvips]{arximspdf}

\doi{10.1214/10-AOP592}
\volume{39}
\issue{4}
\pubyear{2011}
\firstpage{1528}
\lastpage{1543}

\makeatletter

\newtheorem{theorem}{Theorem}[section]
\newtheorem{prop}[theorem]{Proposition}
\newtheorem{cor}[theorem]{Corollary}
\newtheorem{lem}[theorem]{Lemma}

\newproclaim{defn}[theorem]{Definition}
\newproclaim{rmk}[theorem]{Remark}
\newproclaim{ex}[theorem]{Example}

\newcommand{\E}{\mathbf{E}}
\newcommand{\Var}{\operatorname{Var}}

\newcommand{\R}{\mathbb{R}}
\newcommand{\PP}{\mathbf{P}}
\newcommand{\h}{{\widetilde h}}

\newcommand{\iint}{\int\!\!\int}

\makeatother

\begin{document}
\begin{frontmatter}

\title{Concentration of the information in data with log-concave distributions}
\runtitle{Concentration of information}

\begin{aug}
\author[A]{\fnms{Sergey} \snm{Bobkov}\thanksref{t1}\ead[label=e1]{bobkov@math.umn.edu}} and
\author[B]{\fnms{Mokshay} \snm{Madiman}\corref{}\thanksref{t2}\ead[label=e2]{mokshay.madiman@yale.edu}\ead[label=u1,url]{http://www.stat.yale.edu/\textasciitilde mm888}}
\runauthor{S. Bobkov and M. Madiman}
\affiliation{University of Minnesota and Yale University}
\address[A]{School of Mathematics\\
University of Minnesota\\
Minneapolis, Minnesota 55455\\
USA\\
\printead{e1}}
\address[B]{Department of Statistics\\
Yale University\\
New Haven, Connecticut 06511\\
USA\\
\printead{e2}\\
\printead{u1}}
\end{aug}

\thankstext{t1}{Supported in part by NSF Grant DMS-07-06866.}
\thankstext{t2}{Supported in part by a Junior Faculty Fellowship from
Yale University in Spring 2009 and by NSF CAREER Grant DMS-10-56996.}

\received{\smonth{3} \syear{2010}}

%
\begin{abstract}
A concentration property of the functional ${-}\log f(X)$ is demonstrated,
when a random vector $X$ has a log-concave density $f$ on $\R^n$.
This concentration property implies in particular an extension of the
Shannon--McMillan--Breiman strong ergodic theorem to the class of
discrete-time stochastic processes with log-concave marginals.
\end{abstract}

%

\end{frontmatter}

\section{Introduction}

Let $(\Omega,\mathcal{B},\PP)$ be a probability space and let
$X=(X_1,\ldots,\break X_n)$ be a random vector defined on it with each $X_i$
taking values in $\R$. Suppose that the joint distribution of $X$ has
a density $f$ with respect to a reference measure $\nu(dx)$ on $\R^n$.
For most of this paper (except for the purposes of discussion in this
section), the reference measure is simply Lebesgue measure $dx$ on $\R^n$.
The random variable
\[
\h(X) = - \log f(X)
\]
may be thought of as the (random) information content of $X$. Such an
interpretation is well-justified in the discrete case, when $\nu$ is
the counting measure on some countable subset of $\R^n$ on which the
distribution of $X$ is supported. In this case, $\h(X)$ is essentially
the number of bits needed to represent $X$ by a coding scheme that
minimizes average code length~\cite{Sha48}. In the continuous case
(with reference measure $dx$), one may still call $\h(X)$ the information
content even though the coding interpretation no longer holds.
In statistics, one may think of the information content as the
log likelihood function.

The average value of the information content of $X$ is known more commonly
as the entropy. Indeed, the entropy of $X$ is defined by
\[
h(X) = - \int f(x) \log f(x) \,dx = - \E\log f(X).
\]
Observe that we adopt here the usual abuse of notation: we write $h(X)$
even though the entropy is a functional depending only on the distribution
of $X$ and not on the value of $X$.
In general, $h(X)$ may or may not exist (in the Lebesgue sense);
if it does, it takes values in the extended real line $[-\infty
,+\infty]$.

Because of the relevance of the information content in various areas such
as information theory, probability and statistics, it is intrinsically
interesting to understand its behavior. In particular, a natural question
arises: is it true that the information content concentrates around
the entropy in high dimension? In general, there is no reason for such
a concentration property to hold. A main purpose in this note is, however,
to show that when the probability measure on $\R^n$ of interest is
absolutely continuous and log-concave, $\log f(X)$ does possess
a powerful concentration property. Specifically, we prove the following theorem.
\begin{theorem}\label{thm:exp}
Suppose $X=(X_1,\ldots,X_n)$ is distributed according to a log-concave density
$f$ on $\R^n$. Then, for all $t>0$,
\[
\PP\bigl\{ |\h(X) - h(X)| \geq t \sqrt{n} \bigr\} \leq
2 e^{-ct},
\]
where $c>0$ is a universal constant. In fact, one may take $c=1/16$.
\end{theorem}

Note that under the assumption of log-concavity and absolute continuity,
$h(X)$~always exists and is finite (see, e.g.,~\cite{BM09maxent}).

Let us emphasize that the distribution of the difference $\h(X) - h(X)$
is stable under all affine transformations of the space, that is,
\[
\h(TX) - h(TX) = \h(X) - h(X)
\]
for all invertible affine maps $T\dvtx\R^n \rightarrow\R^n$.
In particular, the variance of the information content
\[
\E|\h(X) - h(X)|^2
\]
represents an affine invariant. By Theorem~\ref{thm:exp}, when $f$ is
log-concave,
this variance is bounded by $Cn$ with some universal constant $C$.

In fact, the deviation inequality in Theorem~\ref{thm:exp} amounts to
a stronger bound
$\|\h(X) - h(X)\|_{\psi_1} \leq C\sqrt{n}$ with respect to the
Orlicz norm,
generated by the Young function $\psi_1(t) = e^{|t|} -1$.
This is consistent with the observation that in many standard examples
$\h(X)$ behaves like the sum\vadjust{\goodbreak} of $n$ independent random variables.
For example, when $X$ is standard normal, we have
\[
\h(X) - h(X) = \sum_{i=1}^n \frac{X_i^2 - 1}{2}.
\]
More generally, if $X = (X_1,\ldots,X_n)$ has independent components, then
\[
\h(X) - h(X) = \sum_{i=1}^n \h(X_i) - h(X_i).
\]
These examples show that $\sqrt{n}$-normalization in Theorem \ref
{thm:exp} is chosen
correctly and cannot be improved for the class of log-concave
distributions.

When the dimension $n$ is large, the exponential decay in Theorem \ref
{thm:exp} may be
improved to the Gaussian decay on the interval $0 < t < O(\sqrt{n})$.
\begin{theorem}\label{thm:gaus}
Given a random vector $X$ in $\R^n$ with log-concave density~$f$,
\[
\PP\biggl\{\frac{1}{\sqrt{n}} |{\log f(X)} - \E\log f(X)| \geq
t\biggr\}
\leq3 e^{-ct^2},\qquad 0 \leq t \leq2\sqrt{n} ,
\]
where $c>0$ is a universal constant. In fact, one may take $c=1/16$.
\end{theorem}

Substituting $t = s\sqrt{n}$, rewrite the above inequality as
%
%
\begin{equation}\label{eq:gausconc}
\PP\biggl\{ \biggl|\frac{1}{n} \log
\frac{1}{f(X)} - \frac{h(X)}{n}\biggr| \geq s\biggr\} \leq3 e^{-s^2 n/16},
\end{equation}
for $0 \leq s \leq2$. Equivalently, in terms of the entropy power
$N(X) = \exp\{-\frac{2}{n}\times \E\log f(X)\}$,
we get for the value, say, $s=1$,
\[
\PP\{N(X) e^{-2/n} < f(x)^{2/n} < N(X) e^{2/n}\}
\geq1 - 3 e^{-n/16}.
\]
Thus, with high probability, $f(x)^{2/n}$ is very close to $N(X)$, and
the distribution of $X$ itself is effectively the uniform distribution
on the class of typical observables, or the ``typical set''
[defined to be the collection of all points $x\in\R^n$ such that $f(x)$
lies between $e^{-h(X)-n\varepsilon}$ and $e^{-h(X)+n\varepsilon}$,
for some small fixed $\varepsilon>0$].


A similar concentration inequality was obtained by Klartag and Milman~\cite{KM05},
who compared the value $f(X)$ to the maximum $M$
of the density $f$ and proved that
\[
P\{f(X)^{1/n} > c_0 M^{1/n}\} > 1 - c_1^n
\]
with some absolute constants $c_0, c_1 \in (0,1)$.
Note this result readily follows from Theorem~\ref{thm:gaus},
but not conversely.

Theorems~\ref{thm:exp} and~\ref{thm:gaus}, by entailing an effective uniformity
of the distribution of $X$ on some compact set, provide
a strong, quantitative formulation of the
\textit{asymptotic equipartition property}
for log-concave measures.
To describe this\vadjust{\goodbreak} interpretation, suppose $\mathbb{X}=(X_1,X_2,\ldots)$
is a stochastic process
on the probability space $(\Omega,\mathcal{B},\PP)$, with each $X_i$
taking values in $\R$, and define the corresponding projections
$X^{(n)} = (X_1,\ldots,X_n)$. If $\mathbb{X}$ is stationary, the limit
\[
h(\mathbb{X}) = \lim_{n\rightarrow\infty} \frac{h(X^{(n)})}{n}
\]
exists as long as the increments $h(X^{(n+1)})-h(X^{(n)})$ are finite,
and is called the entropy rate of $\mathbb{X}$. 
For stationary processes $\mathbb{X}$, the question of whether the
information content per coordinate $\frac{\h(X^{(n)})}{n}$ converges
to the limit $h(\mathbb{X})$ (in $L^p$ or in probability or almost surely)
has been extensively studied. In the discrete case, the affirmative answer
to this question goes back to Shannon~\cite{Sha48}, McMillan~\cite{McM53}
and Breiman~\cite{Bre57full}, and the eponymous theorem has been called
``the basic theorem of information theory.''
The continuous case was partially developed by Moy~\cite{Moy611},
Perez~\cite{Per64} and Kieffer~\cite{Kie74}. The definitive version
[almost sure convergence for stochastic processes defined on a standard
Borel space, and allowing more general reference measures $\nu(dx)$ than
Lebesgue and counting measure] is due independently to Barron~\cite{Bar85}
and Orey~\cite{Ore85}; the former in particular gives a clear exposition
and recounting of the history. Specifically, these works imply that if
$\mathbb{X}$ is stationary and ergodic, then, as $n\rightarrow\infty$,
%
%
\begin{equation}\label{eq:aep}
-\frac{1}{n} \log f\bigl(X^{(n)}\bigr) \rightarrow h(\mathbb{X})
\qquad\mbox{a.s.}
\end{equation}
An elementary proof of this fact, called by McMillan the ``asymptotic
equipartition property''
was later given by Algoet and Cover~\cite{AC88}.
For nonstationary processes with arbitrary dependence, the entropy rate
$h(\mathbb{X})$ typically does not exist; so there is no question of
a statement like (\ref{eq:aep}) holding. Nonetheless, together
with Borel--Cantelli's lemma Theorem~\ref{thm:exp} immediately yields
the following
extension 
of the Shannon--McMillan--Breiman phenomenon.
\begin{cor}\label{cor:aep}
Suppose that $\mathbb{X}$ has a log-concave
distribution on $\R^\infty$ with absolutely continuous finite-dimensional
projections. If the limit $h(\mathbb{X})$ exists, the property
(\ref{eq:aep}) holds.
\end{cor}

Note that log-concavity of a probability measure is defined
on arbitrary locally convex spaces via a Brunn--Minkowski type inequality
and is equivalent to the log-concavity of densities of finite-dimensional
projections (in case they are absolutely continuous with respect
to Lebesgue measure; see~\cite{Bor74} for a general theory).
Corollary~\ref{cor:aep} trivially extends to processes
$\mathbb{X}=(X_1,X_2,\ldots)$ where each $X_i$ takes values
in $\R^k$ instead of $\R$, as long as the finite-dimensional projections
$X^{(n)}$ have log-concave distributions.
This, for instance, means that Corollary~\ref{cor:aep} can be applied to
nonstationary Markov chains in $\R^k$ that preserve log-concavity of
the joint distribution
and also have a unique invariant probability measure (the latter
condition ensures
existence of the entropy rate, which can also be easily computed as the mean
under the invariant measure\vadjust{\goodbreak} of the entropy of the conditional density
of $X_2$ given $X_1$).
Furthermore, if the process mixes well enough so that $h(X^{(n)})/n$
converges rapidly to $h(\mathbb{X})$, then Theorem~\ref{thm:gaus} may
be used
to give a convergence rate in probability.

It should also be mentioned that, for Gaussian distributions, tight
concentration inequalities may be derived by simple explicit calculation.
This was done by Cover and Pombra~\cite{CP89} as an ingredient in studying
the feedback capacity of time-varying additive Gaussian noise channels.


The paper is organized in the following way.
As a first step, we consider a one-dimensional version of
Theorem~\ref{thm:exp} (Section~\ref{sec:1d}). In Section \ref
{sec:revlyap}, we
recall some previous work on reverse Lyapunov inequalities, and
present a new variant. It is applied to establish
a concentration property of the logarithm function under what we call
log-concave measures of order $p$ (Sections~\ref{sec:p} and~\ref{sec:log}).
Section~\ref{sec:ls} uses the localization lemma of Lov\'asz and
Simonovits to reduce
the general case to a specific one-dimensional statement.
Section~\ref{sec:pf} completes the proof.


\section{\texorpdfstring{One-dimensional case in Theorem \protect\ref{thm:exp}}{One-dimensional case in Theorem 1.1}}
\label{sec:1d}

We begin by proving the one-dimensional case of Theorem~\ref{thm:exp}.
\begin{prop}\label{prop:1d}
If a random variable $X$ has a log-concave
density $f$, then,
\[
\E e^{({1}/{2}) |{\log f(X)} - \E\log f(X)|} < 4.
\]
\end{prop}
\begin{pf}
Let $X$ be a random variable with log-concave density $f(x)$.
The distribution of $X$ is supported on some interval
$(a,b)$, finite or not, where $f$ is positive and $\log f$ is concave.
Introduce the function
\[
I(t) = f(F^{-1}(t)),\qquad 0<t<1,
\]
where $F^{-1}\dvtx(0,1) \rightarrow(a,b)$ is the inverse to the distribution
function $F(x) = \PP\{X \leq x\}$, $a < x < b$. The function $I$ is positive
and concave on $(0,1)$ and uniquely determines $F$ up to a shift
parameter (\cite{Bob961}, Proposition A.1).

Given a function $\Psi= \Psi(u,v)$, write a general identity
\[
\iint\Psi(f(x),f(y)) f(x) f(y) \,dx \,dy =
\int_0^1 \int_0^1 \Psi(I(t),I(s)) \,dt \,ds.
\]
In particular, for any $\alpha\in[0,1)$,
%
%
\begin{equation}\label{eq:subst}\qquad
\iint e^{\alpha|{\log f(x)} - \log f(y)|} \,dF(x) \,dF(y)
=
\int_0^1 \int_0^1 e^{\alpha|{\log I(t)} - \log I(s)|} \,dt \,ds.
\end{equation}
Here the right-hand side does not change when multiplying $I$ by
a positive scalar, so one may assume that $I(1/2) = 1/2$. But then,
by concavity of $I$, we have
\[
\min\{t,1-t\} \leq I(t) \leq1.
\]
From this,
\begin{eqnarray*}
\log I(t) - \log I(s) &\leq&- \log\min\{s,1-s\}, \\
\log I(s) - \log I(t) &\leq&- \log\min\{t,1-t\},
\end{eqnarray*}
so
\[
|{\log I(t)} - \log I(s)| \leq- \log\min\{t,s,1-t,1-s\}.
\]

Hence, the right-hand side of (\ref{eq:subst}) does not exceed
\begin{eqnarray*}
\int_0^1 \int_0^1 e^{-\alpha\log\min\{t,s,1-t,1-s\}} \,dt \,ds
&=& 4 \int_0^{1/2} \int_0^{1/2} \min\{t,s\}^{-\alpha} \,dt \,ds \\
&=& \frac{2^{1 + \alpha}}{(1 - \alpha) (2 - \alpha)}.
\end{eqnarray*}

Finally, by Jensen's inequality with respect to $dF(y)$, the left-hand side
of (\ref{eq:subst}) majorizes
\[
\int e^{\alpha|{\log f(x)} - \int\log f(y) \,dF(y)|} \,dF(x) =
\E e^{\alpha|{\log f(X)} - \E\log f(X)|} ,
\]
so that we have
%
%
\begin{equation}\label{eq:gen-alpha}
\E e^{\alpha|{\log f(X)} - \E\log f(X)|} \leq\frac{2^{1 + \alpha
}}{(1 - \alpha) (2 - \alpha)}.
\end{equation}
Choosing the value $\alpha= 1/2$, and
observing that $\frac{8}{3} \sqrt{2} < 4$,
we may conclude.
\end{pf}

Note also that a direct application of Chebyshev's inequality yields
\[
\PP\{ |{\log f(X)} - \E\log f(X)| \geq t\} \leq4 e^{-t/2}
\]
for all $t>0$. While the exponent here is slightly better than that in
Theorem~\ref{thm:exp}, we make no effort here (or anywhere in this
paper) to come
up with optimal constants.


\section{Reverse Lyapunov inequalities}
\label{sec:revlyap}

Given a random variable $\eta> 0$, the Lyapunov inequality states that
%
%
\begin{equation}\label{eq:lyap}
\lambda_a^{b-c} \lambda_c^{a-b} \geq\lambda_b^{a-c},\qquad
a \geq b \geq c > 0,
\end{equation}
where $\lambda_p = \E\eta^p$ is the moment function of $\eta$. Equivalently,
it expresses a well-known and obvious property that the function
$p \rightarrow\log\lambda_p$ is convex on the positive half-axis $p>0$.

What is less obvious, for certain classes of probability distributions
on $(0,+\infty)$, the inequality (\ref{eq:lyap}) may be reversed
after a suitable
normalization of the moment function. In particular, when $\eta$
has a distribution with increasing hazard rate (in particular, if $\eta$
has a log-concave density), then as was shown by
Barlow, Marshall and Proschan (\cite{BMP63}, page 384), we have
%
%
\begin{equation}\label{eq:bmp-lyap}
\bar{\lambda}_a^{b-c} \bar{\lambda}_c^{a-b} \leq\bar{\lambda
}_b^{a-c},\qquad
a \geq b \geq c \geq1 \qquad(c \mbox{ is integer}),
\end{equation}
for the normalized moment function
\[
\bar{\lambda}_p = \frac{1}{\Gamma(p+1)} \E\eta^p.
\]
Note that $\bar{\lambda}_p = 1$ for all $p>0$ for the standard exponential
distribution, which thus plays an extremal role in this class.

This result has many interesting applications. For example, applying
it to the parameters $a=p+1$, $b = p$, $c = p-1$, we have
%
%
\begin{equation}\label{eq:bmp-lyap2}
\E\eta^{p+1} \cdot\E\eta^{p-1} \leq\biggl(1 + \frac{1}{p}
\biggr) (\E\eta^p)^2,
\end{equation}
provided that $p \geq2$ is integer. If the distribution of $\eta$ is
log-concave, the case $p=1$ can also be included in this inequality,
which is due to a Khinchine-type inequality by Karlin, Proschan and
Barlow~\cite{KPB61}, namely,
\[
\E\eta^p \leq\Gamma(p+1) (\E\eta)^p,\qquad p \geq1.
\]

However, in some problems, it is desirable to remove the requirement
that $c$ is integer in (\ref{eq:bmp-lyap}). This is implied by results
of Borell~\cite{Bor73a}
for the class of log-concave densities. To be more precise, he proved
the following (Theorem~2 in~\cite{Bor73a}).
%
%
\begin{prop}\label{prop31}
Let $\eta$ be a nonnegative concave function,
defined on an open convex body $\Omega\subset\R^n$. Then the function
\[
p \longrightarrow\frac{(p+1) \cdots(p+n)}{n!} \int_\Omega\eta
(x)^p \,dx
\]
is log-concave in $p \geq0$.
\end{prop}

To relate this to (\ref{eq:bmp-lyap}), let us start with a continuous
convex function
$u\dvtx\Delta\rightarrow\R$, defined on some closed segment
$\Delta\subset(0,+\infty)$, such that $e^{-u(x)}$ is a probability
density. For large $n$, consider convex bodies
\[
\Omega_n = \biggl\{(x_1,\ldots,x_n,x) \in\R_+^n \times\Delta\dvtx
x_1 + \cdots+ x_n \leq1 - \frac{u(x)}{n}\biggr\}.
\]
Their volumes satisfy, as $n \rightarrow\infty$,
%
%
\begin{equation}\label{eq:revlyap-mid1}
n! |\Omega_n| =
\int_\Delta\biggl(1 - \frac{u(x)}{n}\biggr)^n \,dx
\rightarrow\int_\Delta e^{-u(x)} \,dx = 1,
\end{equation}
and for every $p \geq0$,
%
%
\begin{equation}\label{eq:revlyap-mid2}\qquad
v_n(p) =
\frac{1}{|\Omega_n|} \int_{\Omega_n} x^p \,dx_1 \cdots dx_n \,dx
\quad\rightarrow\quad
v(p) =
\int_\Delta x^p e^{-u(x)} \,dx.\vadjust{\goodbreak}
\end{equation}
By Proposition~\ref{prop31}, applied to $\eta(x_1,\ldots,x_n,x) = x$,
the functions
\[
w_n(p) = \frac{(p+1) \cdots(p+n)}{n^{p+1} n!} v_n(p),\qquad p
\geq0,
\]
are log-concave, so the limit will also be a log-concave
function, if it exists. (Note that we have added a log-linear factor
$n^{p+1}$.) But
\[
\frac{(p+1) \cdots(p+n)}{n^{p+1} n!} \rightarrow\frac{1}{\Gamma(p+1)}.
\]
Therefore, in view of (\ref{eq:revlyap-mid1}) and (\ref
{eq:revlyap-mid2}), the resulting limit
$\frac{1}{\Gamma(p+1)} v(p)$ represents a log-concave function, as well.

On this step, the assumption that $u$ was defined on a closed segment
can be relaxed, and we arrive at the following corollary
(which seems not to be mentioned in~\cite{Bor73a} or anywhere else).
\begin{cor}\label{cor:rev-lyap}
If a random variable $\eta> 0$ has a log-concave
distribution, then the function
\[
\bar{\lambda}_p = \frac{1}{\Gamma(p+1)} \E\eta^p,\qquad p \geq0,
\]
is log-concave. Equivalently, we have a reverse Lyapunov's inequality
%
%
\begin{equation}\label{eq:rev}
\bar{\lambda}_a^{b-c} \bar{\lambda}_c^{a-b} \leq\bar{\lambda
}_b^{a-c},\qquad
a \geq b \geq c \geq0.
\end{equation}
\end{cor}

In connection with the concentration problem and the Kannan--Lov\'
asz--Simonovits conjecture
within the class of spherically symmetric distributions on~$\R^n$,
reverse Lyapunov's inequalities were considered in~\cite{Bob03gafa1}.
The following alternative variant of Corollary~\ref{cor:rev-lyap} is
proposed there.
\begin{prop}\label{prop:rev-lyap}
Given a random variable $\eta> 0$ with a
log-concave distribution, the
function $\hat{\lambda}_p = \E(\frac{\eta}{p})^p$ is log-concave
in $p > 0$,
and therefore satisfies (\ref{eq:rev}).
\end{prop}

This is proved in~\cite{Bob03gafa1} by an application of the Pr\'
ekopa--Leindler inequality,
and is perhaps more convenient for applications involving asymptotics.

There is much more that can be (and has been) said about reverse Lyapunov
inequalities; a gentle introduction may be found in~\cite{BM10norm}.


\section{Log-concave distributions of order $p$}
\label{sec:p}

\begin{defn}
A random variable $\xi>0$ will be said to have a
log-concave distribution of order $p \geq1$, if it has a density of
the form
\[
f(x) = x^{p-1} g(x),\qquad x>0,
\]
where the function $g$ is log-concave on $(0,+\infty)$.
\end{defn}

When $p=1$, we obtain the class of all (nondegenerate) log-concave
probability distributions on $(0,+\infty)$.

The meaning of the parameter $p$ is that it is responsible for
a strengthened concentration. For example, the inequality (\ref
{eq:bmp-lyap2}), which
holds by Corollary~\ref{cor:rev-lyap} for all real $p \geq1$, may
equivalently be
rewritten in terms of $\xi$ as
%
%
\begin{equation}\label{eq:var1}
\Var(\xi) \leq\frac{1}{p} (\E\xi)^2.
\end{equation}
Alternatively, if we start with Proposition~\ref{prop:rev-lyap} and
apply (\ref{eq:rev})
with $a=p+1$, $b = p$, $c = p-1$ ($p > 1$), we get
$\E\eta^{p+1} \E\eta^{p-1} \leq C_p (\E\eta^p)^2$ with constants
$C_p = (p+1)^{p+1} (p-1)^{p-1} p^{-2p}$. Equivalently,
%
%
\begin{equation}\label{eq:var2}
\Var(\xi) \leq(C_p - 1) (\E\xi)^2
\end{equation}
in the class of log-concave $\xi$ of order $p$.
Asymptotically $C_p = 1 + \frac{1}{p} + O(\frac{1}{p^3})$, as
$p \rightarrow+\infty$, so the bound (\ref{eq:var2}) is very close
to (\ref{eq:var1})
for large values of $p$.
\begin{ex}
Let $\xi$ have a Gamma distribution with shape parameter $p$ 
(where $p>0$ is real), that is, with density
\[
f(x) = \frac{1}{\Gamma(p)} x^{p-1} e^{-x},\qquad x>0.
\]
It is log-concave if and only if $p \geq1$, in
which case $p$ will be the order of log-concavity for this distribution.
Note that $\E\xi= \Var(\xi) = p$, and (\ref{eq:var1}) becomes
equality. Hence,
the factor $1/p$ in (\ref{eq:var1}) is optimal.
\end{ex}
\begin{prop}\label{prop:var-p}
If $\xi> 0$ has a log-concave distribution of
order $p \geq1$, then
\[
\Var(\log\xi) \leq\frac{d^2}{dp^2} \log\Gamma(p).
\]
Equality is attained at the Gamma distribution with shape parameter
$p$. 
\end{prop}
\begin{pf}
Write the density of $\xi$ as $f(x) = x^{p-1} g(x)$
with log-concave $g$. One may assume that $g$ is a density, as well.
Indeed, otherwise consider random variables $\xi_c = c \xi$ $(c>0)$.
Then $\Var(\log\xi_c) = \Var(\log\xi)$ and $\xi_c$ has density
\[
f_c(x) = c^{-p} x^{p-1} g(x/c) = x^{p-1} g_c(x),
\]
where $g_c(x) = c^{-p} g(x/c)$. Since $f$ decays at infinity exponentially
fast, the same is true for $g$. Hence, $g$ is integrable, and one can
choose $c$ such that $\int g_c(x) \,dx = 1$. So the reduction to the case
where $g$ is a density is achieved.

Thus, let $g$ be a log-concave probability density, such that
$f(x) = x^{p-1} g(x)$ is the density of $\xi$. Consider a random variable\vadjust{\goodbreak}
$\eta> 0$ with density $g$. Then, by Corollary~\ref{cor:rev-lyap},
the function
\[
u(q) = \log\E\eta^{q-1} - \log\Gamma(q),\qquad q \geq0,
\]
is concave. Differentiating twice with respect to $q$, we get
\[
u''(q) = \frac{\E\eta^{q-1} \log^2 \eta-
(\E\eta^{q-1} \log\eta)^2}{(\E\eta^{q-1})^2} -
\frac{d^2}{dq^2} \log\Gamma(q) \leq0.
\]
But at the point $q = p$, we have
\[
\E\eta^{p-1} = \int x^{p-1} g(x) \,dx = \int f(x) \,dx = 1,
\]
and so
\begin{eqnarray*}
u''(p) + \frac{d^2}{dp^2} \log\Gamma(p) &
= &
\E\eta^{q-1} \log^2 \eta- (\E\eta^{q-1} \log\eta)^2
\\
&
= &
\int x^{p-1} \log^2 x g(x) \,dx -
\biggl(\int x^{p-1} \log x g(x) \,dx\biggr)^2 \\
&=& \Var(\log\xi).
\end{eqnarray*}
Proposition~\ref{prop:var-p} is proved.
\end{pf}

It is to be noted that the right-hand side in Proposition~\ref{prop:var-p}
is the trigamma function, which has the alternate representation
\[
\psi_1(p)=\sum_{n=1}^{\infty} \frac{1}{(n+p)^2} ,
\]
and behaves like $1/p$ for large values of $p$. Hence,
%
%
\begin{equation}\label{eq:lcp-conc}
\Var(\log\xi) \leq\frac{C}{p}
\end{equation}
with some absolute constant $C$ (in fact, one may take $C=1$).
%
This can also be seen by using
Proposition~\ref{prop:rev-lyap}. Indeed, the same argument as above yields
%
%
\begin{equation}\label{eq:lcp-conc2}
\Var(\log\xi) \leq\frac{d^2}{dp^2} (p-1) \log(p-1) = \frac{1}{p-1},
\end{equation}
which holds for any $p>1$. Here the right-hand side has an incorrect
behavior when $p$ is close to 1. In fact, for all log-concave $\xi$,
we have
%
%
\begin{equation}\label{eq:lc-conc}
\Var(\log\xi) \leq C
\end{equation}
with some absolute constant $C$. For the proof, one can apply, for example,
Borell's concentration lemma (\cite{Bor74}, Lemma 3.1).
Together with (\ref{eq:lc-conc}), (\ref{eq:lcp-conc2}) also yields
the bound (\ref{eq:lcp-conc}).

In the proof of Theorem~\ref{thm:exp}, we use the values $p=n$, the dimension
of the space. Since the one-dimensional case can be treated separately
(rather easily), the assumption $p \geq2$ can be made in applications.
\begin{rmk}\label{rmk:aff}
The notion of a log-concave measure of order $p$ may be
extended in a natural way to the class of one-dimensional log-concave
probability measures $\mu$ on $\R^n$. More precisely, we say that
$\mu$
has order $p$, if $\mu$ is supported on some interval $\Delta\subset
\R^n$,
bounded or not, and has a density there of the form
\[
\frac{d\mu(x)}{dx} = \ell(x)^{p-1} g(x),\qquad x \in\Delta,
\]
where $\ell$ is a positive affine function on $\Delta$,
$g$ is log-concave on $\Delta$, and where $dx$ stands for the Lebesgue
measure on this interval.
In this case, the inequality (\ref{eq:lcp-conc}) and other similar
results should
be properly read in terms of $\ell$. For example, we have
$\Var(\log\ell) \leq\frac{C}{p}$ with respect to $\mu$.
\end{rmk}


\section{Concentration of the logarithm function}
\label{sec:log}

It is natural to try to sharpen Proposition~\ref{prop:var-p} and the
resulting asymptotic
bound (\ref{eq:lcp-conc}) in terms of deviations of $\log\xi$ from
its mean or quantiles.

Let $\xi> 0$ be a random variable with log-concave distribution
of order $p+1$, that is, with density of the form
\[
f(x) = x^p g(x),\qquad x>0,
\]
where $p \geq0$ and $g$ is a log-concave function. Let $\zeta$ be an
independent copy of $\xi$. Then for all $\alpha\in[0,p]$,
\begin{eqnarray*}
\E e^{\alpha| {\log\xi}- \log\zeta|} & = &
2 \E e^{\alpha(\log\xi- \log\zeta)} 1_{\{\xi> \zeta\}} \\
& \leq&
2 \E e^{\alpha(\log\xi- \log\zeta)}
=
2 \E\xi^\alpha\E\zeta^{-\alpha}\\
&=&
2 \int x^{p + \alpha} g(x) \,dx \int x^{p - \alpha} g(x) \,dx.
\end{eqnarray*}

The quantity $\E e^{\alpha| {\log\xi}- \log\zeta|}$ does not change
if we multiply $\xi$ and $\zeta$ by a positive scalar. Hence, as in the
proof of Proposition~\ref{prop:var-p}, we may assume that $g$ is a
probability density
of some random variable, say, $\eta$. Applying Jensen's inequality, we thus
conclude that
%
%
\begin{equation}\label{eq:log-mid}
\E e^{\alpha| {\log\xi}- \E\log\xi|} \leq
2 \E\eta^{p + \alpha} \E\eta^{p - \alpha},\qquad
0 \leq\alpha\leq p,
\end{equation}
provided that $\E\eta^p = 1$ (which means that $f$ is a density).
But by the reverse Lyapunov's inequality of Corollary \ref
{cor:rev-lyap}, applied
with $a = p + \alpha$, $b = p$, $c = p - \alpha$, we obtain that
\[
\E\eta^{p + \alpha} \E\eta^{p - \alpha} \leq
\frac{\Gamma(p + \alpha+ 1) \Gamma(p - \alpha+ 1)}{\Gamma(p + 1)^2}.
\]
Note that when $\alpha= 1$, this inequality returns us to inequality
(\ref{eq:var1}).

Thus, from (\ref{eq:log-mid}),
\[
\E e^{\alpha| {\log\xi}- \E\log\xi|} \leq2
\frac{\Gamma(p + \alpha+ 1) \Gamma(p - \alpha+ 1)}{\Gamma(p +
1)^2},\qquad
0 \leq\alpha\leq p.\vadjust{\goodbreak}
\]

The right-hand side here seems perhaps not quite convenient to deal with,
especially when $p \pm\alpha$ are not integer.
Alternatively, it might be better to use Proposition \ref
{prop:rev-lyap}, which gives
\[
\E e^{\alpha| {\log\xi}- \E\log\xi|} \leq2
\frac{(p + \alpha)^{p + \alpha} (p - \alpha)^{p -
\alpha}}{p^{2p}},\qquad
0 \leq\alpha\leq p.
\]
Indeed, write
\[
\frac{(p + \alpha)^{p + \alpha} (p - \alpha)^{p - \alpha
}}{p^{2p}} =
\biggl(1 - \frac{\alpha^2}{p^2}\biggr)^{p - \alpha}
\biggl(1 + \frac{\alpha}{p}\biggr)^{2\alpha}.
\]
The first factor on the right may be bounded just by 1. For the second one,
using $(1 + t)^{1/t} \leq e$ ($t \geq0$), one has
\[
\biggl(1 + \frac{\alpha}{p}\biggr)^{2\alpha} =
\biggl(1 + \frac{\alpha}{p}\biggr)^{({p}/{\alpha})
({2\alpha^2}/{p})}
\leq e^{2\alpha^2/p}.
\]
Therefore, we have a preliminary Gaussian estimate:
\[
\E e^{\alpha|{\log\xi}- \E\log\xi|} \leq2 e^{2\alpha^2/p},\qquad
0 \leq\alpha\leq p .
\]

Similarly, one may also obtain a one-sided estimate, since
like in inequality (\ref{eq:log-mid}) we also have
\[
\E e^{\alpha(\log\xi- \log\zeta)} =
\E\eta^{p + \alpha} \E\eta^{p - \alpha},\qquad
0 \leq|\alpha| \leq p,
\]
provided that $\E\eta^p = 1$. These estimates are collected below
after replacing $p$ by $p-1$ for convenience.
\begin{lem}\label{lem:logconc}
If $\xi> 0$ has a log-concave distribution of
order $p >1$, then
%
%
\begin{eqnarray}\label{eq:log1}
\E e^{\alpha|{\log\xi}- \E\log\xi|} &\leq&2 e^{{2\alpha
^2}/({p-1})},\qquad
0 \leq\alpha\leq p-1,
\\
%
%
\label{eq:log2}
\E e^{\alpha(\log\xi- \E\log\xi)} &\leq&
e^{{2\alpha^2}/({p-1})},\qquad
0 \leq|\alpha| \leq p-1.
\end{eqnarray}
\end{lem}

In particular, we obtain for log-concave densities of order $p$ on the
positive half-line
a $p$-dependent version of Proposition~\ref{prop:1d} (which was stated
for log-concave densities
on the line).
\begin{cor}\label{cor:log}
If $\xi> 0$ has a log-concave distribution of
order $p \geq1$, then
\[
\E e^{({1}/{6}) \sqrt{p} |{\log\xi}- \E\log\xi|} < 3 .
\]
\end{cor}
\begin{pf}
First, assume $p \geq2$ and choose $\alpha= c\sqrt{p}$ in (\ref
{eq:log1}) with
$0 < c \leq1/\sqrt{2}$ (so that $\alpha\leq p-1$). Then, using
$p/(p-1) \leq2$, we have
\[
\E e^{c\sqrt{p} |{\log\xi}- \E\log\xi|} \leq2 e^{4c^2}.
\]
Taking, for example, $c = 1/6$, the right-hand side will not exceed
$2 e^{1/9} < 3$. Hence,
\[
\E e^{({1}/{6}) \sqrt{p} |{\log\xi}- \E\log\xi|} < 3.\vadjust{\goodbreak}
\]

For the remaining range $1 \leq p < 2$,
one has $\sqrt{p}/6 < 1/4$, and we have by Proposition~\ref{prop:1d}
[or more precisely, inequality (\ref{eq:gen-alpha})] that
%
%
\begin{equation}
\E e^{({\sqrt{p}}/{6}) |{\log\xi}- E \log\xi|} <
\E e^{({1}/{4}) |{\log\xi}- E \log\xi|} \leq\frac{2^{5/4}}{3/4
\times7/4} < 2.
\end{equation}
Thus, the desired statement is proved with a uniform bound of 3.
\end{pf}

Observe that Proposition~\ref{prop:1d} corresponds to $p=1$, and that
while it clearly applies
as stated to log-concave densities of order $p$ (since these are
subclasses of the log-concave densities),
Corollary~\ref{cor:log} with the additional $\sqrt{p}$ term in the
exponent provides the correct generalization for large $p$.



\section{Reduction to dimension one}
\label{sec:ls}

To reduce Theorems~\ref{thm:exp} and~\ref{thm:gaus} to a specific
statement about dimension one
(in fact---about log-concave distributions of order $p = n$), we
apply a localization argument of Lov\'asz and Simonovits~\cite{LS93}.
More precisely, we need one variant of the localization lemma, proposed
in~\cite{KLS95}, Corollary 2.4, which we state with minor modification
as a lemma.
\begin{lem}\label{lem:ls}
Let $g$ and $h$ be integrable continuous functions on
a bounded open convex set $\Omega$ in $\R^n$, such that
\[
\int_\Omega g(x) \,dx > 0,\qquad \int_\Omega h(x) \,dx = 0.
\]
Then for some interval $\Delta\subset\Omega$ and a positive
affine function $\ell$ on $\Delta$,
\[
\int_\Delta g \ell^{n-1} > 0,\qquad \int_\Delta h \ell^{n-1} = 0,
\]
where the integrals are with respect to Lebesgue measure on $\Delta$.
\end{lem}

Equivalently, given that $\int_\Omega h(x) \,dx = 0$, if
for all couples $(\Delta,\ell)$ with $\int_\Delta h\times \ell^{n-1} = 0$,
we have that
\[
\int_\Delta g \ell^{n-1} \leq0,
\]
then
\[
\int_\Omega g(x) \,dx \leq0.
\]
This formulation enables the desired-dimensional reduction.
\begin{lem}\label{lem:ls-red}
Suppose $X$ is a random vector taking values in an open convex set
$\Omega$ in $\R^n$,
where it has a positive continuous density $f$, such that $\E|{\log f(X)}|$ is finite.
Let $\mu_\ell$ denote a probability measure on a line segment $\Delta
\subset\Omega$ with density
\[
f_\ell(x) = \frac{1}{Z} f(x) \ell(x)^{n-1},\vadjust{\goodbreak}
\]
where $\ell$ is a positive affine function, defined on $\Delta$,
and $Z = \int_\Delta f(x) \ell(x)^{n-1} \,dx$ is a normalizing constant.
Given $\alpha\geq0$ and $A \geq1$, if for any such one-dimensional
measure $\mu_\ell$, we have
%
%
\begin{equation}\label{eq:if-1}
\E_{ \ell} e^{({\alpha}/{\sqrt{n}}) |{\log f} - \E_\ell\log
f|} \leq A,
\end{equation}
where $\E_{ \ell}$ stands for the expectation with respect to $\mu
_\ell$,
then
%
%
\begin{equation}\label{eq:then-n}
\E\exp\biggl\{\frac{\alpha}{\sqrt{n}} |{\log f(X)} - \E\log
f(X)|\biggr\}
\leq A .
\end{equation}
%
\end{lem}
\begin{pf}
Without loss of generality, take $\Omega$ to be bounded,
and assume that
$\E\log f(X) = 0$, or in other words,
%
%
\begin{equation}\label{eq:ent0}
\int_\Omega\log f(x) f(x) \,dx = 0.
\end{equation}
In this case, (\ref{eq:then-n}) becomes
%
%
\begin{equation}\label{eq:then-ag}
\int_\Omega\bigl(e^{({\alpha}/{\sqrt{n}})
|{\log f(x)}|} - A\bigr) f(x) \,dx \leq0.
\end{equation}
This corresponds to Lemma~\ref{lem:ls} with
\[
h(x) = \log f(x) f(x) \quad\mbox{and}\quad
g(x) = \bigl(e^{({\alpha}/{\sqrt{n}}) |{\log f(x)}|} - A\bigr) f(x).
\]
Hence, to derive (\ref{eq:then-ag}) under (\ref{eq:ent0}), it
suffices to take an arbitrary
interval $\Delta\subset\Omega$ and a positive affine function $\ell$
on $\Delta$, such that
%
%
\begin{equation}\label{eq:alm1}
\int_\Delta\log f(x) f(x) \ell(x)^{n-1} \,dx = 0,
\end{equation}
and to show that
%
%
\begin{equation}\label{eq:alm2}
\int_\Delta\bigl(e^{({\alpha}/{\sqrt{n}})
|{\log f(x)}|} - A\bigr) f(x) \ell(x)^{n-1} \,dx \leq0.
\end{equation}

Using the definition of $\mu_\ell$, inequalities (\ref
{eq:alm1}) and (\ref{eq:alm2}) take the form
\[
\int\log f d\mu_\ell= 0,
\int\bigl(e^{({\alpha}/{\sqrt{n}}) |{\log f}|} - A\bigr)
\,d\mu_\ell\leq0,
\]
which can be written together as (\ref{eq:if-1}).
\end{pf}


\section{\texorpdfstring{Proofs of Theorems \protect\ref{thm:exp} and \protect\ref{thm:gaus}}
{Proofs of Theorems 1.1 and 1.2}}
\label{sec:pf}

Keeping the same notation as in the previous section, first note that
\[
\log f - \E_{ \ell} \log f = (\log f_\ell- \E_{ \ell} \log
f_\ell) -
(n-1) (\log\ell- \E_{ \ell} \log\ell),
\]
so
\[
|{\log f} - \E_{ \ell} \log f| \leq
|{\log f_\ell}- \E_{ \ell} \log f_\ell| +
(n-1) |{\log\ell}- \E_{ \ell} \log\ell|.
\]
By convexity of the functional $\xi\rightarrow\log\E e^\xi$, we
have that
%
%
\begin{eqnarray}\label{eq:collect}\quad
\log\E_{ \ell} e^{({\alpha}/({2\sqrt{n}})) |{\log f} - \E_\ell
\log f|}
&\leq&
\tfrac{1}{2}
\log\E_{ \ell} e^{({\alpha}/{\sqrt{n}}) |{\log f_\ell}- \E
_\ell\log f_\ell|}\nonumber\\[-8pt]\\[-8pt]
&&{}+
\tfrac{1}{2}
\log\E_{ \ell} e^{({\alpha(n-1)}/{\sqrt{n}}) |{\log\ell}
- \E_\ell\log\ell|}.\nonumber
\end{eqnarray}

Since $f_\ell$ is the density of the one-dimensional log-concave probability
measure~$\mu_\ell$, by Proposition~\ref{prop:1d}, whenever
$0 \leq\alpha\leq\frac{1}{2} \sqrt{n}$,
%
%
\begin{equation}\label{eq:piece1}
\E_{ \ell} e^{({\alpha}/{\sqrt{n}}) |{\log f_\ell}- \E_\ell
\log f_\ell|}
< 4.
\end{equation}

To estimate the second expectation in (\ref{eq:collect}), it is useful
to note that $\mu_\ell$
has order $p = n$ (cf. Remark~\ref{rmk:aff}). If $n=1$, this
expectation is just 1.
If $n \geq2$, by the inequality (\ref{eq:log1}) of Lemma~\ref{lem:logconc},
we have
%
%
\begin{equation}\label{eq:piece2}
\E_{ \ell} e^{({\alpha(n-1)}/{\sqrt{n}}) |{\log\ell}- \E
_\ell\log\ell|}
\leq
2 e^{2\alpha^2 (n-1)/n} \leq2 e^{2\alpha^2},
\end{equation}
provided that $0 \leq\alpha\leq\sqrt{n}$. This bound automatically
holds for $n=1$, as well.

Collecting the bounds (\ref{eq:piece1}) and (\ref{eq:piece2}) in (\ref
{eq:collect}), we get that,
for all $0 \leq\alpha\leq\frac{1}{2} \sqrt{n}$,
\[
\log\E_{ \ell} e^{({\alpha}/({2\sqrt{n}})) |{\log f} - \E_\ell
\log f|}
\leq
\tfrac{1}{2} \log(8 e^{2\alpha^2}).
\]
Hence, using $\sqrt{8} < 3$ (to simplify the constant),
\[
\E_{ \ell} e^{({\alpha}/({2\sqrt{n}})) |{\log f} - \E_\ell\log
f|} \leq
3 e^{\alpha^2}.
\]
Now, replace $\alpha$ with $2\alpha$. We then get that
\[
\E_{ \ell} e^{({\alpha}/{\sqrt{n}}) |{\log f} - \E_\ell\log
f|} \leq
3 e^{4\alpha^2}, \qquad 0 \leq\alpha\leq\tfrac{1}{4} \sqrt{n}.
\]
Recalling Lemma~\ref{lem:ls-red} (whose assumptions hold
for all log-concave densities), we arrive at the following theorem.
\begin{theorem}\label{thm:mgf}
Given a random vector $X$ in $\R^n$ with
log-concave density $f(x)$,
\[
\E\exp\biggl\{\frac{\alpha}{\sqrt{n}} |{\log f(X)} - \E\log
f(X)|\biggr\}
\leq3 e^{4\alpha^2},\qquad 0 \leq\alpha\leq\frac{1}{4} \sqrt{n}.
\]
\end{theorem}

Choose $\alpha= 1/4$. Denoting
$\xi= \frac{1}{4\sqrt{n}} |{\log f(X)} - \E\log f(X)|$,
we have
$\E e^\xi\leq3 e^{1/4}$. Hence,
$\E e^{\xi/4} \leq3^{1/4} e^{1/16} < 2$. This gives the following.
\begin{cor}\label{cor:mgf}
Given a random vector $X$ in $\R^n$ with
log-concave density $f(x)$,
\[
\E\exp\biggl\{\frac{1}{16\sqrt{n}} |{\log f(X)} - \E\log
f(X)|\biggr\}
\leq2.
\]
\end{cor}

By applying Chebyshev's inequality, we arrive at Theorem~\ref{thm:exp}
with $c=1/16$.
From Theorem~\ref{thm:mgf}, by Chebyshev's inequality, we also have
\[
\PP\biggl\{\frac{1}{\sqrt{n}} |{\log f(X)} - \E\log f(X)| \geq
t\biggr\}
\leq3 e^{4\alpha^2 - \alpha t},\qquad t>0,
\]
provided that $0 \leq\alpha\leq\frac{1}{4} \sqrt{n}$. Taking the optimal
value $\alpha= t/8$ gives Theorem~\ref{thm:gaus}.

%

%
\printaddresses

\end{document}